\documentclass[11pt]{amsart}

\usepackage{amsfonts, amstext, amsmath, amsthm, amscd, amssymb}
\usepackage{graphicx, epsfig}

\let\oldmarginpar\marginpar
\renewcommand\marginpar[1]{\oldmarginpar[\raggedleft\footnotesize #1]
{\raggedright\footnotesize #1}}

\newcommand{\Z}{\mathbb{Z}}
\newcommand{\s}{\sigma}

\renewcommand{\c}{\cite}

\newcommand{\kb}[1]{\ensuremath{\langle #1 \rangle}}
\newcommand{\asplice}{\ensuremath{\makebox[0.3cm][c]{\raisebox{-0.3ex}{\rotatebox{90}{$\asymp$}}}}}
\newcommand{\del}{\partial}
\newcommand{\G}{\mathbb{G}}
\renewcommand{\H}{\mathbb{H}}
\newcommand{\C}{\mathcal{C}}
\renewcommand{\d}{\delta}
\def \BRT {{Bollob\'as--Riordan--Tutte }}
\newcommand{\Q}{\mathbb{Q}}
\renewcommand{\H}{\mathbb{H}}
\newcommand{\TT}{\mathbb{T}}

\theoremstyle{plain}
\newtheorem{theorem}{Theorem}[section]
\newtheorem{corollary}[theorem]{Corollary}

\newtheorem{prop}[theorem]{Proposition}

\newtheorem{conjecture}[theorem]{Conjecture}

\theoremstyle{definition}

\newtheorem{question}{Question}

\newtheorem*{namedtheorem}{\theoremname}
\newcommand{\theoremname}{testing}

\title[A survey on the Turaev genus of knots]{A survey on the Turaev genus of knots}
\author[A.\ Champanerkar]{Abhijit Champanerkar}
\address{Department of Mathematics, College of Staten Island \& The Graduate Center, City University of New York, New York, NY}
\email{abhijit@math.csi.cuny.edu, ikofman@math.csi.cuny.edu}
\author[I. \ Kofman]{Ilya Kofman}

\begin{document}

\maketitle

\begin{abstract}
The Turaev genus of a knot is a topological measure of how far a given
knot is from being alternating.  Recent work by several authors has
focused attention on this interesting invariant.  We discuss how the
Turaev genus is related to other knot invariants, including the Jones
polynomial, knot homology theories, and ribbon-graph polynomial
invariants.
\end{abstract}

\section{Introduction}

Knots and links have been studied using graphs associated to their
diagrams since the first knot tables were compiled in the late 1800's.
Separately, ribbon graphs, which are cellularly embedded graphs on a
two-dimensional surface, have a long history, not only in graph theory
but in the study of Riemann surfaces, Galois theory, quantum field
theory and many other subjects (see \cite{LZ} for example).  Only
recently, though, have ribbon graphs been associated to knot diagrams
in a way that yields powerful new invariants of knots and links.  In
\cite{DFKLS1}, Dasbach, Futer, Kalfagianni, Lin and Stoltzfus
discovered that the Jones polynomial is a specialization of the \BRT
polynomial of a particular ribbon graph on a surface originally
constructed by Turaev, called the Turaev surface of a knot.  The
minimal genus of such a surface, called the Turaev genus of a knot, is
an interesting new invariant that measures how far a given knot (or
link) is from being alternating.

The aim of this paper is to give some historical background about the
ideas leading to the Turaev surface and Turaev genus, explain the
connections to the more well-known graphs associated to knot diagrams,
review some modern applications in knot homology theories, and lastly
focus on open problems and new research directions related to
the Turaev genus. A natural question is how the Turaev genus is
related to other diagrammatic, geometric and topological invariants of
knots and links.

\subsection{Definition}
Let $D$ be a diagram of a link $L$.  For any crossing
\includegraphics[height=0.3cm, angle=90]{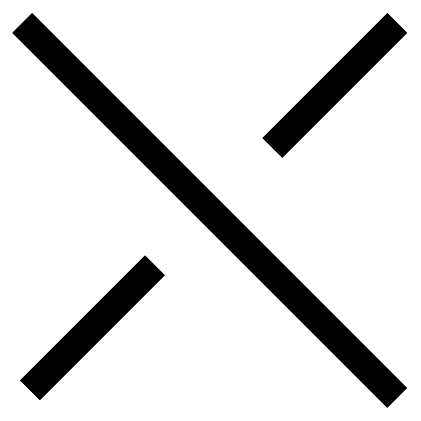} in $D$, we
obtain the $A$--smoothing as $\asplice$ and the $B$--smoothing as
$\asymp$.  A state $s$ of $D$ is a choice of smoothing at every
crossing, resulting in a disjoint union of circles in the plane. Let
$|s|$ denote the number of circles in $s$.  Let $s_A$ denote the
all--$A$ state, for which every crossing of $D$ has an $A$--smoothing.
Similarly, $s_B$ is the all--$B$ state of $D$.

\begin{center}
\begin{figure}[h]
\includegraphics[height=1in]{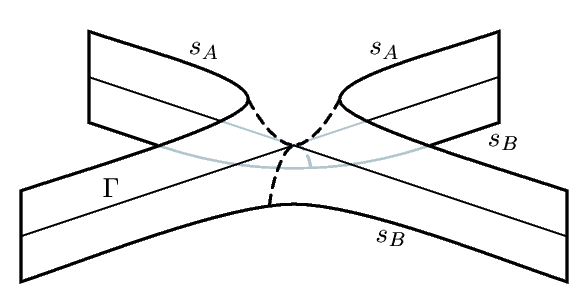}\ \ \ \ 
\includegraphics[height=1.2in]{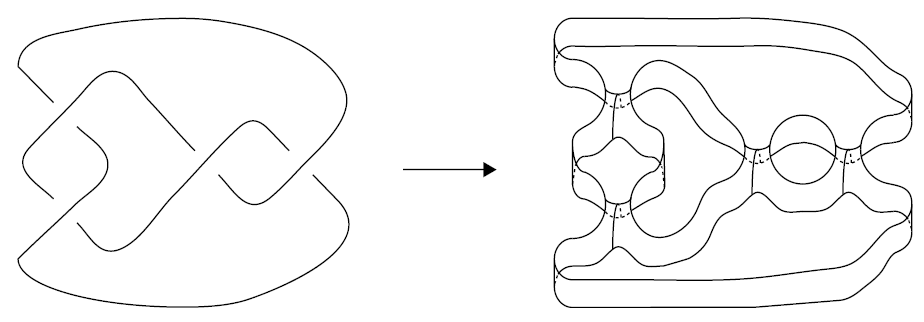}
\caption{Cobordism between $s_A$ and $s_B$ (figures from \c{cromwell} and \c{Abe1})}
\label{figure:turaev-surface}
\end{figure}
\end{center}

Now, at every crossing of $D$, we put a saddle surface which bounds
the $A$--smoothing on the top and the $B$--smoothing on the bottom as
shown in Figure \ref{figure:turaev-surface}.  In this way, we get a
cobordism between $s_A$ and $s_B$, with the link projection at the
level of the saddles. See Figure~\ref{figure:turaev-surface}.

For any diagram $D$, the {\em Turaev surface} $F(D)$ is obtained by
attaching $|s_A|+|s_B|$ discs to all boundary circles. 
% The description above gives a Morse decomposition of $F(D)$, and 
The {\em Turaev genus} of $D$ is defined by
$$g_T(D) = g(F(D)) = (c(D) + 2 - |s_A| - |s_B|)/2.$$
The {\em Turaev genus} of any non-split link $L$ is defined by
$$ g_T(L) = \min\{ g_T(D)\ | \ D\  {\rm is\ a\ diagram\ of\ }L \}. $$

The properties below follow easily from the definitions. See
\cite{DFKLS1} for proofs and figures.
\begin{enumerate}
\item[(a)] $F(D)$ is an unknotted closed orientable surface in $S^3$; i.e., $S^3-F(D)$ is a disjoint union of two handlebodies.
\item[(b)] $D$ is alternating on $F(D)$.
\item[(c)] $L$ is alternating if and only if $g_T(L)=0$, and if $D$ is an alternating diagram then $F(D) = S^2$.
\item[(d)] The projection of $D$ is a 4-valent graph which gives a cell decomposition of $F(D)$, for which the 2-cells can be
  checkboard colored on $F(D)$, with discs corresponding to $s_A$ and
  $s_B$ respectively colored white and black.
\item[(e)] This cell decomposition is a Morse decomposition of $F(D)$, for
  which $D$ and the crossing saddles are at height zero, and the $s_A$ and
  $s_B$ 2-cells are the maxima and minima, respectively.
\end{enumerate}

For example, any non-alternating pretzel knot  can be
made alternating on the torus as follows:

\centerline{\includegraphics[height=1.5in]{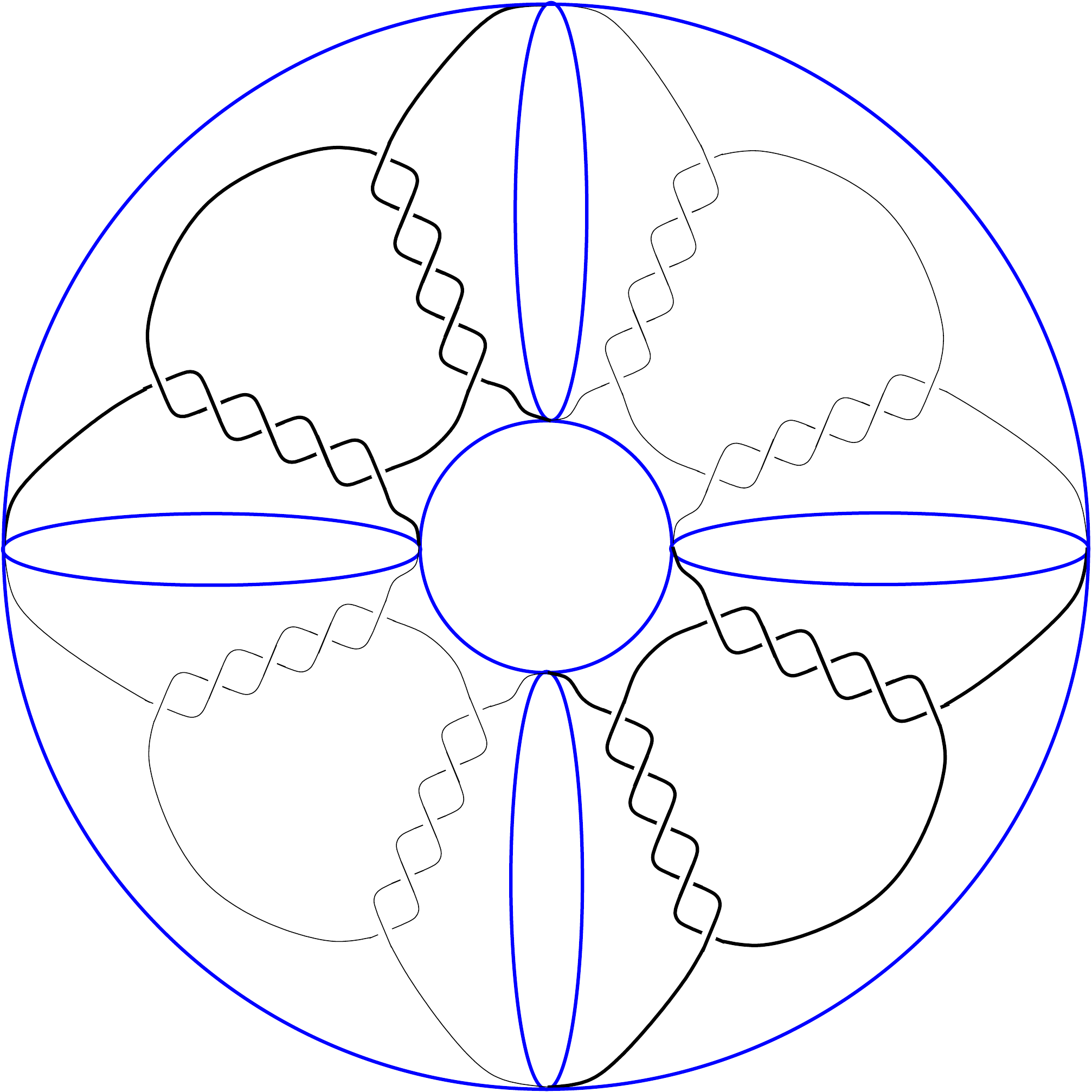}}

In fact, the Turaev genus of any non-alternating pretzel knot (or more
generally, any non-alternating Montesinos knot) equals one, and its
Turaev surface is the torus.  
We will return to this fact in Section \ref{bounds}. 
See \cite{hajij} for a nice animation of Turaev surfaces. 

The paper is organized as follows: we discuss the motivation behind
the Turaev surface in Section 2.  In Section 3, we discuss the Turaev
surface in the context of ribbon graphs.  In Section 4, we discuss
applications to knot homology theories.  In Section 5, we discuss
known bounds for the Turaev genus.  In Section 6, we discuss research
directions and open problems related to the Turaev genus.

\section{Tait's Conjecture}

Modern knot theory began in late 1800's when Tait, Little and others
tried to make a periodic table of elements by tabulating knot diagrams
by crossing number.  The only invariants at this time were of the
form, ``minimize something among all diagrams,'' such as crossing
number, unknotting number, bridge number, etc.  Such invariants are
easy to define but hard to compute: Diagrams that are minimal with
respect to one property may not be minimal with respect to other
properties.

\begin{center}
\begin{figure}[h]

\includegraphics[height=0.9in]{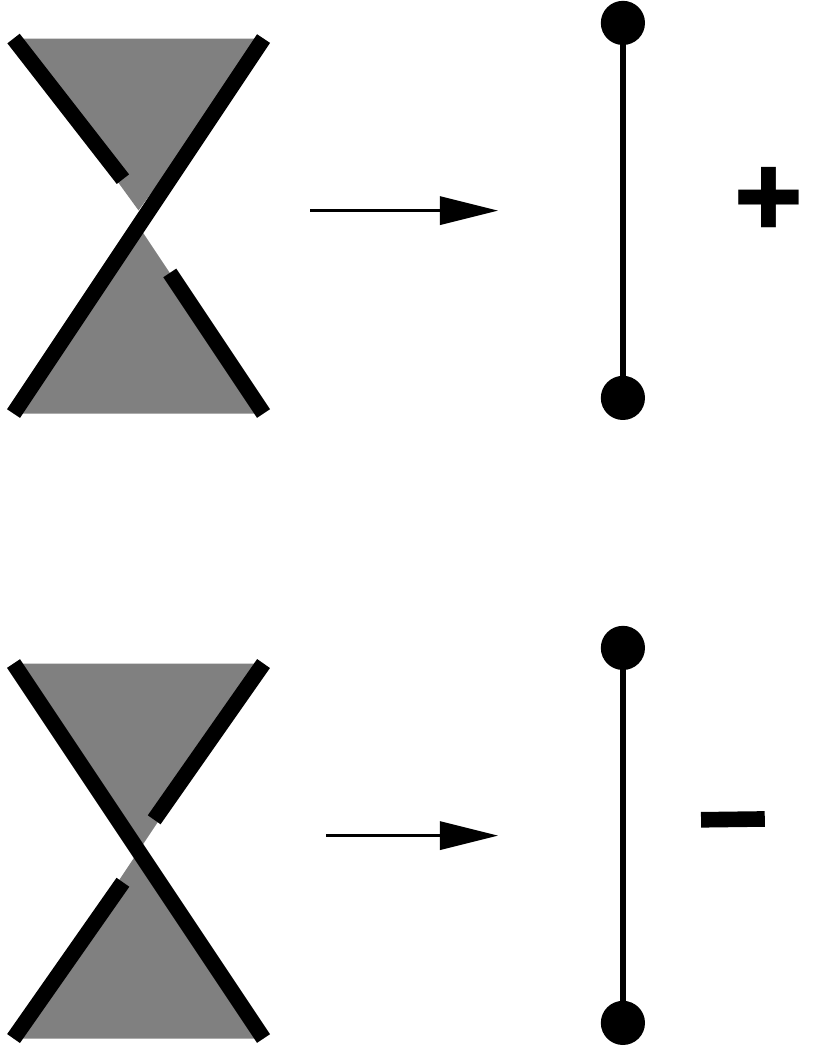}
\hspace*{1.2cm}
\includegraphics[height=0.9in]{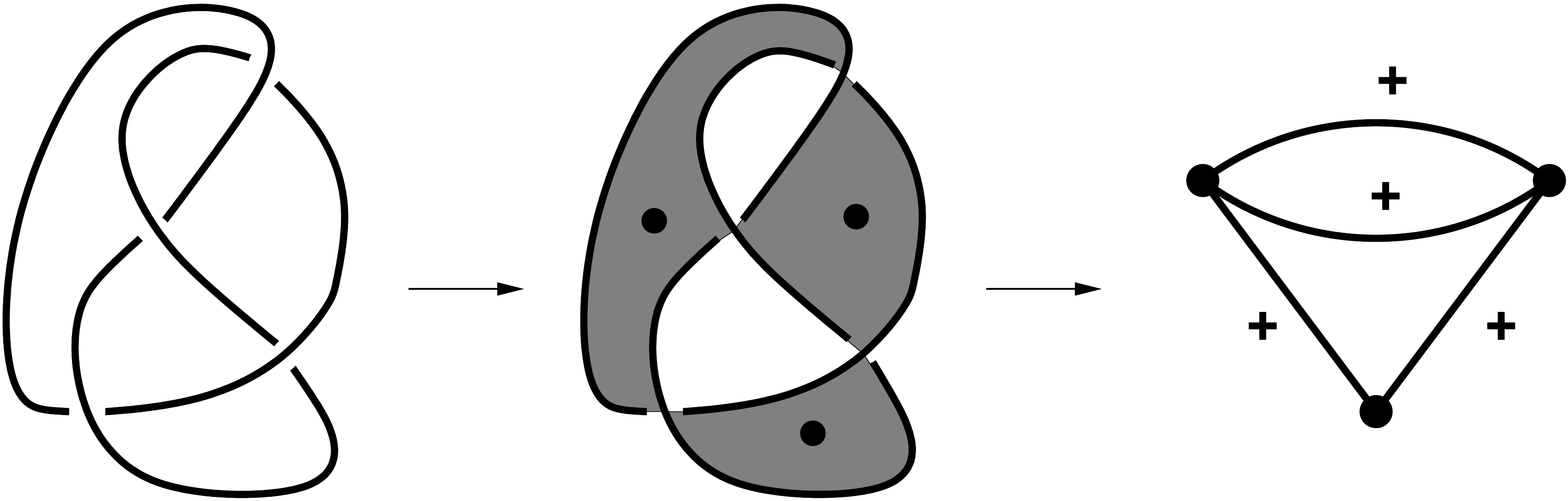}\\
\ \\
\includegraphics[height=1in]{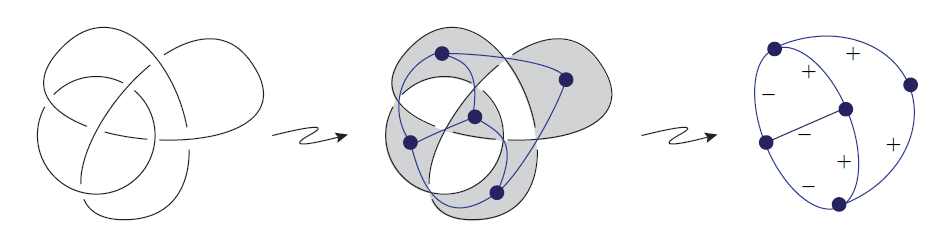}
\label{figure:Tait graph}
\caption{Edge sign convention and Tait graphs (lower figure from \c{Moffatt_duals})}
\end{figure}
\end{center}

There is a correspondence between connected link diagrams $D$
and connected plane graphs $G$ with signed edges.  It follows from the
Jordan Curve Theorem that any link diagram can be checkerboard
colored.  The {\em Tait graph} $G$ of $D$ is obtained by checkerboard
coloring complementary regions of $D$, assigning a vertex to every
shaded region, an edge to every crossing and a $\pm$ sign to every
edge as in Figure \ref{figure:Tait graph}.

Conversely, we can recover the diagram from any signed 
planar graph by taking its medial graph, and making 
crossings according to the sign on each edge: \\
\centerline{\includegraphics[height=1cm]{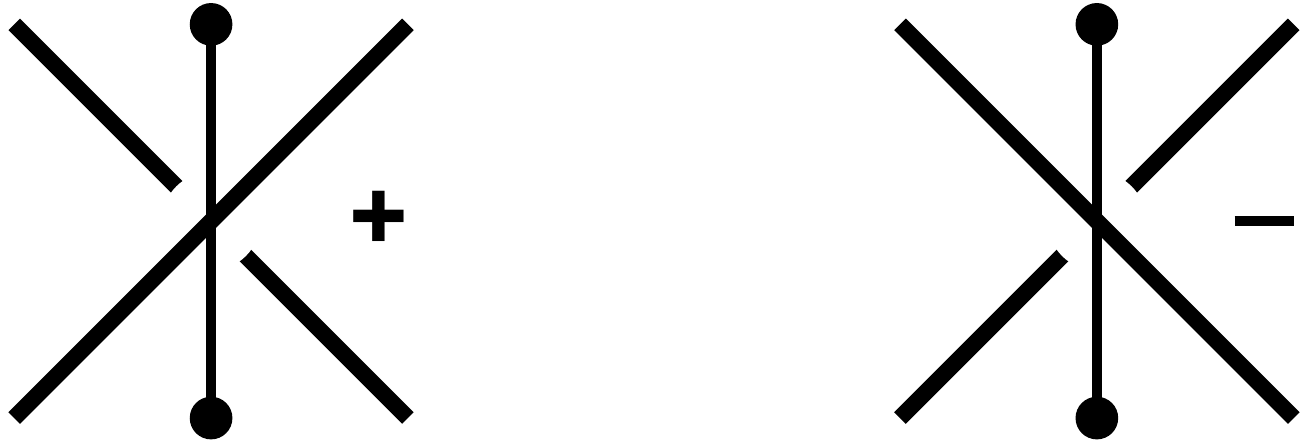}} \\
Tait graphs for opposite checkerboard colorings are planar duals.

A link diagram is {\em alternating} if the crossings alternate between
overcrossing and undercrossing as one walks along every component of
the link. A link is {\em alternating} if it has a reduced alternating
diagram. Tait emphasized the importance of an alternating diagram, for
which all the edges in its Tait graph have the same sign, so it
corresponds to an unsigned plane graph, determined by the diagram up
to planar duality.

\begin{conjecture}[Tait]
An alternating link always has an alternating diagram that has minimal
crossing number among all diagrams for that link.
\end{conjecture}

A proof had to wait about 100 years until the Jones polynomial (1984),
which led to several new ideas that were used to prove Tait's
conjecture \cite{kauffman, murasugi, Thistlethwaite}.  Below we
follow a later proof by Turaev \cite{turaev} using Turaev surfaces
defined above.  In Section \ref{KHsection}, we discuss the spanning
tree expansion for the Jones polynomial introduced in
\cite{Thistlethwaite}, which is of independent interest.

The simplest combinatorial approach to the Jones polynomial is via the
Kauffman bracket $\kb{D} \in \Z[A,A^{-1}]$ defined recursively by
\begin{enumerate}
\item $\kb{\includegraphics[height=0.3cm, angle=90]{crossing}} = A \; \kb{\asplice} + A^{-1}  \kb{\asymp} $
\item $\kb{\bigcirc\ D} = \d \; \kb{D},\quad \d=-A^{-2}-A^2 $ 
\item $\kb{\bigcirc} = 1 $
\end{enumerate}

For any link $L$ with a diagram $D$ of writhe $w(D)$, the Jones
polynomial is determined by the Kauffman bracket as
$V_L(A^{-4})=(-A)^{-3w(D)}\kb{D}$. We will use $\kb{L}$ or $V_L(t)$
depending on the context.

Besides this axiomatic definition, Kauffman \cite{kauffman} expressed
$\kb{L}$ as a sum over all possible states of $L$: If $L$ has $n$
crossings, all possible $A$ and $B$ smoothings yield $2^n$ states $s$.
Let $a(s)$ and $b(s)$ be the number of $A$ and $B$ smoothings,
respectively, to get $s$.  
$$ \kb{L} = \sum_{{\rm states}\ s} A^{a(s)-b(s)}\ (-A^2 - A^{-2})^{|s|-1} $$

A diagram $D$ is {\em adequate} if ($i$) $|s_A|>|s|$ for any state
$s$ with exactly one $B$--smoothing, and ($ii$) $|s_B|>|s|$ for any
state $s$ with exactly one $A$--smoothing.  In particular, any reduced
alternating diagram is adequate.  A link is adequate if it has an
adequate diagram.

The proof of Tait's Conjecture now follows from three claims (see \cite{cromwell}):  
\begin{enumerate}
\item[($i$)] Although defined for diagrams, the Jones polynomial $V_L(t)$ is a
link invariant.
\item[($ii$)] $s_A$ and $s_B$ contribute the extreme terms $\pm t^\alpha$ and $\pm t^\beta$ of $V_L(t)$, 
which determine the span $V_L(t) = \alpha - \beta$, which is a link invariant.
In particular,
$$ \max\deg_A\kb{D} - \min\deg_A\kb{D} \leq 2(c(D)+|s_A(D)|+|s_B(D)|-2) $$
with equality if $D$ is adequate, hence if $D$ is alternating. 
\item[($iii$)] By Turaev's {\em dual-state lemma}, $|s_A(D)|+|s_B(D)| = 2+c(D)-2g_T(D)$.  Thus,
$$ \max\deg_A\kb{D} - \min\deg_A\kb{D} \leq 4c(D) - 4g_T(D) $$
%% It follows that 
\begin{equation}\label{spaneq} 
{\rm span}\, V_L(t) \leq c(L) - g_T(L) 
\end{equation}
with equality if $L$ is adequate, hence if $L$ is alternating.  
If $D$ is a prime non-alternating diagram, then $g_T(D)>0$ so we get a strict inequality. %: 
%% $\max\deg_A\kb{D} - \min\deg_A\kb{D} < 4c(D)$.
Thus, span $V_L(t)= c(L)$ if and only if $L$ is alternating, from which Tait's Conjecture follows.
\end{enumerate}
%
%% We define span $V_L(t) = \alpha - \beta$, which is a link invariant.  It follows that 
%% \begin{equation}\label{spaneq} 
%% {\rm span}\, V_L(t) \leq c(L) - g_T(L) 
%% \end{equation}
%% with equality if $L$ is adequate.  
%% Thus, span $V_L(t)= c(L)$ if and only if $L$ is alternating, from which Tait's Conjecture follows.

Therefore, for any adequate link $L$ with an adequate diagram $D$, we
have (see \cite{Abe1}):
\begin{equation}\label{adeqeq} 
 g_T(L) = g_T(D) = \frac{1}{2}\left(c(D)-|s_A(D)|-|s_B(D)|\right)+1 = c(L) - {\rm span}\, V_L(t)
\end{equation}

\subsection{Spanning trees and Jones polynomial}
The Tait graph plays a role in an earlier proof of the Tait conjecture via
Thistlethwaite's spanning tree expansion of the Jones polynomial.
Thistlethwaite \c{Thistlethwaite} gave an expansion of $V_L(t)$ in
terms of the spanning trees of the Tait graph of any diagram of $L$.
Every spanning tree contributes a monomial to $V_L(t)$. 

For non-alternating diagrams, these monomials may cancel with each
other, but for alternating diagrams, such cancelations do not occur.
Thus, for alternating links, the number of spanning trees is exactly
the $L^1$-norm of coefficients of $V_L(t)$, and the span of $V_L(t)$
is maximal, equal to the crossing number.  This gives a different
proof of claim (3) above for alternating links.  Thistlethwaite also
showed that the Jones polynomial of an alternating link can be
obtained as a specialization of the Tutte polynomial of its Tait
graph.

The Tutte polynomial is a fundamental and ubiquitous invariant of
graphs, which can be defined by a state sum over all subgraphs, by
contraction-deletion operations, and by a spanning tree expansion, any
of which could have led to the Jones polynomial three decades earlier!
Tutte's original definition in \cite{Tutte} used the spanning tree
expansion, which relies on the concept of activity of edges with
respect to a spanning tree.

In Section 4, we discuss the spanning tree expansion and its
applications in more detail.  For example, for the figure-eight knot, Figure \ref{figure:unknots}
shows the skein resolution tree in terms of spanning trees.

\section{Ribbon graphs and polynomial invariants} 

In this section we look at two graph-theoretic generalizations of the
ideas above.  

First, Turaev's construction gives rise to a graph embedded on a
surface, i.e. a ribbon graph, in a way that generalizes the Tait graph
of a diagram $D$. When $D$ is non-alternating, the Tait graph must
have signs to encode all the crossing information of D. Instead, we
can construct an un-signed ribbon graph whose topology completely
encodes $D$. We will formally define ribbon graphs below, which can be
more general graphs on surfaces.

Second, Thistlethwaite's specialization of the Tutte polynomial to the
Jones polynomial of alternating links also generalizes in the ribbon
graph setting: a specialization of the \BRT polynomial gives the Jones
polynomial for any link.  Moreover, the spanning trees of a plane
graph also have natural counterparts in the ribbon graph setting. We
discuss all of these ideas below.

\subsection{Ribbon graphs}

An {\em oriented ribbon graph} is a cellularly embedded graph in an
oriented surface (precisely, a multi-graph for which loops and
multiple edges are allowed) that is embedded in such a way that its
complement is a union of open discs on the surface.  A ribbon graph is
also described as a band decomposition by thickening the cellularly
embedded graph.  See Table \ref{table:ribbon-graphs} and Figure
\ref{figure:ribbon-graph-from-state}.  The embedding, combined with
the orientation on the surface, determines {\em a cyclic order} on the
edges at every vertex, and also a cell structure for the surface.
Terms for the same or closely related objects include: combinatorial
maps, fat graphs, cyclic graphs, graphs with rotation systems, ribbon
and arrow marked graphs and dessins d'enfant (see \c{BR1, MMbook, LZ}
and references therein).

A ribbon graph $\G$ can be considered both as a geometric and as a
combinatorial object.  The combinatorial definition is given as
follows: let $(\s_0,\,\s_1,\,\s_2)$ be permutations of
$\{1,\ldots,2n\}$, such that $\s_1$ is a fixed-point free involution
and $\s_0\,\s_1\,\s_2=1$.  We define the orbits of $\s_0$ to be the
vertex set $V(\G)$, the orbits of $\s_1$ to be the edge set $E(\G)$,
and the orbits of $\s_2$ to be the face set $F(\G)$.  Let $v(\G)$,
$e(\G)$ and $f(\G)$ be the numbers of vertices, edges and faces of
$\G$.  The preceding data determine an embedding of $\G$ on a closed
orientable surface, denoted $S(\G)$, as a cell complex.  The set
$\{1,\ldots,2n\}$ can be identified with the directed edges (or
half-edges) of $\G$.  Thus, $\G$ is connected if and only if the group
generated by $\s_0,\,\s_1,\,\s_2$ acts transitively on
$\{1,\ldots,2n\}$.  The genus of $S(\G)$ is called the genus of $\G$,
$g(\G)$.  If $\G$ has $k(\G)$ components,
$2g(\G)=2k(\G)-v(\G)+e(\G)-f(\G)=k(\G)+n(\G)-f(\G)$, where
$n(\G)=e(\G)-v(\G)+k(\G)$ denotes the nullity of $\G$.  Henceforth, we
assume that $\G$ is a connected, orientable ribbon graph.
See Table \ref{table:ribbon-graphs} for an example of distinct ribbon
graphs with the same underlying graph.

\begin{table}[h]
\begin{center}
\begin{tabular}{ccc}
\includegraphics[width=0.8in]{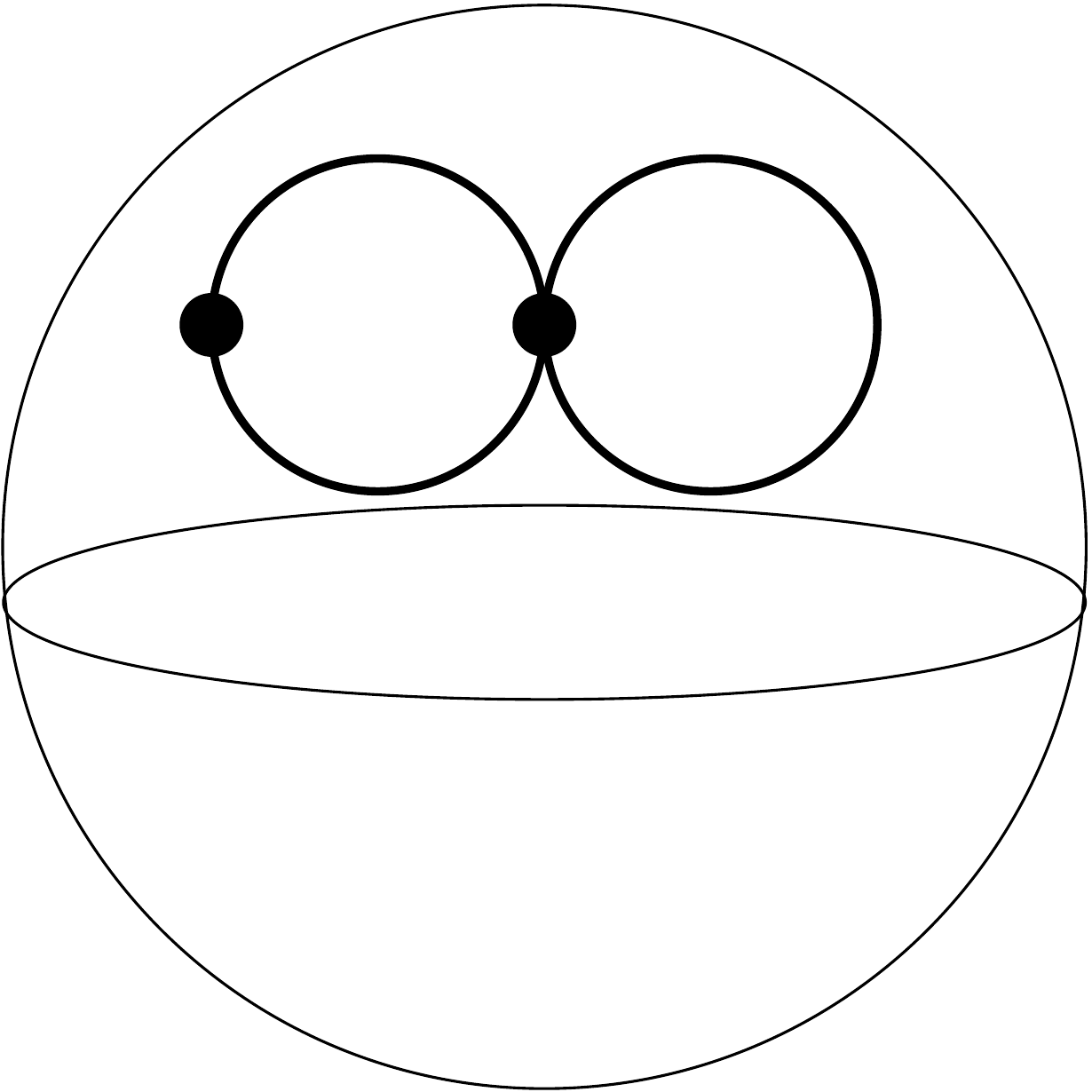} &
\includegraphics[width=1.1in]{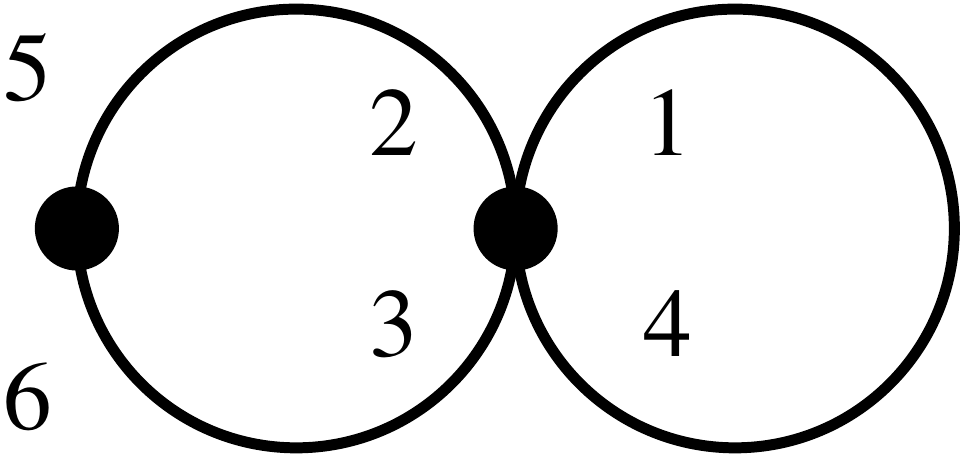} &
\begin{tabular}{l}
 $\sigma_0=(1234)(56)$ \\
 $\sigma_1=(14)(25)(36)$ \\
 $\sigma_2=(246)(35)$ \\
 \\  \\ 
\end{tabular}
\\
\includegraphics[width=1.35in]{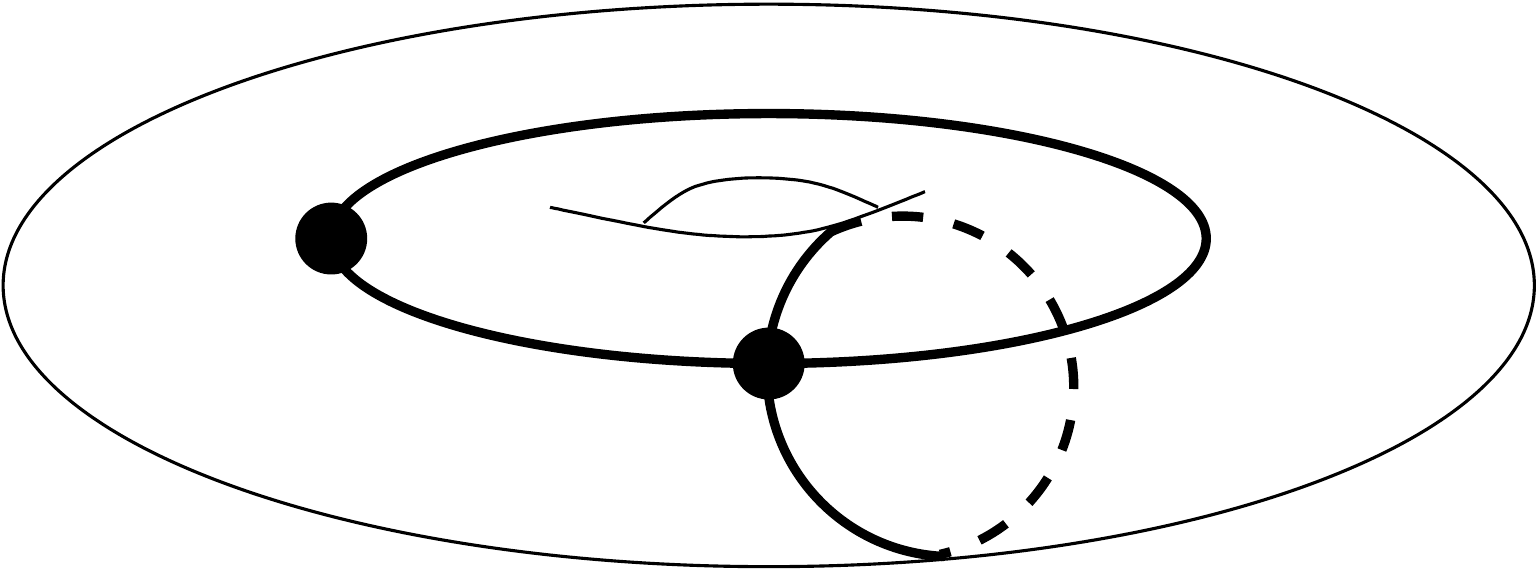} &
\includegraphics[width=1.1in]{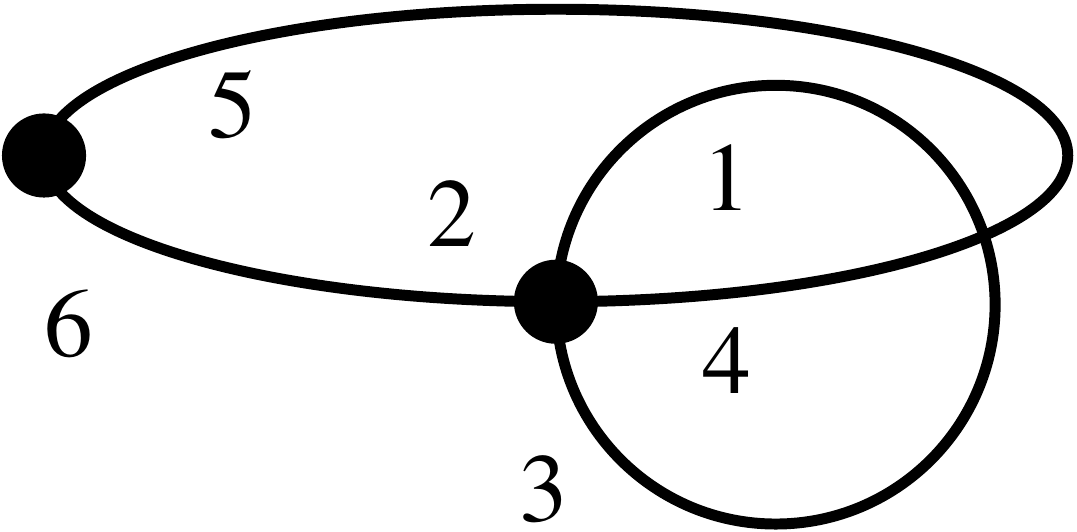} &
\begin{tabular}{l}
$\sigma_0=(1234)(56)$ \\
$\sigma_1=(13)(26)(45)$\\
$\sigma_2=(152364)$ \\
\\ \\ 
\end{tabular}
\end{tabular}
\end{center}
\caption{Ribbon graphs described as graphs on surfaces, as combinatorial maps and as permutations.}
\label{table:ribbon-graphs}
\end{table}

A ribbon graph $\H$ is a {\it ribbon subgraph} of $ \G$ if $\H$ can be
obtained by deleting vertices and edges of $\G$. A ribbon subgraph
$\H\subset\G$ is called a {\it ribbon spanning subgraph} if
$V(\H)=V(\G)$. Note that the surface on which $\H$ is cellularly embedded 
need not be the same surface on which $\G$ is cellularly embedded 
(i.e. $S(\H)$ need not be the same as $S(\G)$), and $g(\H)\leq g(\G)$. 

Bollob\'as and Riordan \c{BR1} extended the Tutte polynomial to a
polynomial invariant of oriented ribbon graphs $C(\G)\in\Z[X,Y,Z]$ in
a way that takes into account the topology of the ribbon graph $\G$.
The Bollob\'as--Riordan--Tutte polynomial has a spanning ribbon subgraph
expansion given by the following sum:
$$ C(\G) = \sum_{\H\subseteq\G} (X-1)^{k(\H)-k(\G)}\; Y^{n(\H)}\;Z^{g(\H)}.$$

In \c{BR2}, they generalized it to a four-variable polynomial
invariant $R(\G)$ of non-orientable ribbon graphs.

\subsection{Ribbon graphs from link diagrams}

Turaev's construction  gives rise to a ribbon graph in a way that
generalizes the Tait graph of a diagram $D$.  The projection of $D$ can be
checkerboard colored on the Turaev surface $F(D)$ with $|s_A|$ white
regions (at height $> 0$), and $|s_B|$ black regions (at height $<
0$). Let $\G_A$ (and similarly $\G_B$) be the graph on $F(D)$ obtained by
assigning a vertex to every white region (respectively, black region),
and an edge to every crossing as we do for the Tait graph. Then the 
complementary regions of $\G_A$ (similarly of $\G_B$) on $F(D)$ are the 
$s_B$ circles (respectively, $s_A$ circles) which bound discs on $F(D)$. 
 Note that
$$ v(\G_A)=|s_A|=f(\G_B),\quad e(\G_A)=e(\G_B)=c(D),\quad 
f(\G_A)= |s_B|=v(\G_B)$$

Thus $\G_A,\,\G_B$ are dual ribbon graphs embedded in $F(D)$ with $
g(\G_A)=g(\G_B)=g_T(D) $.  If $D$ is alternating, $\G_A$ and $\G_B$
are Tait graphs which are planar duals on $F(D)=S^2$. 
If $D$ is $A$--adequate ($B$--adequate), as defined in Section 2, then
$\G_A\ (\G_B)$ has no loops. 

We can also obtain $\G_A$ directly from the link diagram as follows:
 \begin{enumerate}
 \item For a given diagram $D$, use the $A$--smoothing of every crossing to obtain the
state $s_A$, add a ribbon edge (band) joining the two arcs at every smoothed
crossing.
\item Checkerboard color complementary regions of the circles of $s_A$, and orient the circles as the oriented boundary of the black regions.
%% Orient the circles of $s_A$ consistently to agree with the orientation of the plane.
\item Collapse each state circle of $s_A$ to a vertex of $\G_A$,
  preserving the cyclic order of the ribbon edges.
\item Order the half edges at each vertex using the ordering on the
  crossings. The ordering gives us the permutations describing $\G_A$.

 \end{enumerate}
 
\begin{figure}[h]
\begin{center}
  \includegraphics[width=5in]{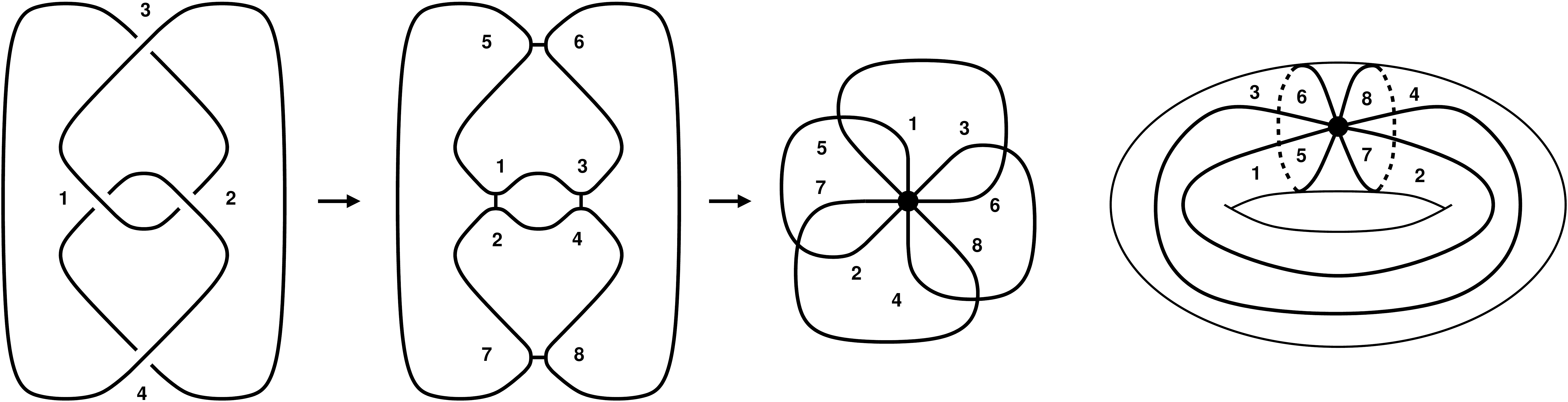}
\end{center}
\caption{Ribbon graph $\G_A$ for a four-crossing diagram of the trefoil knot}
\label{figure:GA}
\end{figure}

This is illustrated in Figure \ref{figure:GA}. $\G_B$ can be obtained similarly 
by starting from the all-B state $s_B$. 

Thistlethwaite \cite{Thistlethwaite} showed that if $L$ is
alternating, then $V_L(t) \doteq T_G(-t,-1/t)$, where $G$ is the Tait
graph of $L$ and $T_G(x,y)$ is its Tutte polynomial.  
In \cite{DFKLS1}, it was shown that the Jones polynomial $V_L(t)$ is a
specialization of the \BRT polynomial of the all--$A$ ribbon graph
$\G_A$.  Chmutov \cite{Chmutov} extended these ideas to virtual links
and non-orientable ribbon graphs, and to links given as a diagram on a
surface.  In \c{MMbook}, a unified description is given for all these
knot and ribbon graph polynomial invariants using the four-variable polynomial $R(\G)$:

\begin{theorem}[\cite{DFKLS1, Chmutov}] Let $D$ be either a classical
  link diagram, a link diagram on a surface, or a virtual link diagram
  and $\G$ be the all--$A$ ribbon graph of $D$. Then
$$ \kb{D} = \delta^{k(\G)-1} A^{n(\G)-r(\G)} R(\G; -A^4, A^{-2}\delta, \delta^{-1}, 1) $$
where $\kb{D}$ is the Kauffmann bracket of $D$ and $\delta = -A^2-A^{-2}$.
\end{theorem}

\subsection{Ribbon graphs for dual states}
States $s$ and $\overline{s}$ are called dual states if the smoothing
at every crossing in $\overline{s}$ is opposite to that in $s$,
e.g. $s_A$ and $s_B$ are dual states.  The Turaev surface construction
applies to any pair of dual states to obtain the surface $F(D_s)$,
with the projection of $D$ emebedded as a 4-valent cellular graph. As
before, we get dual ribbon graphs $G_s$ and $G_{\overline{s}}$ as Tait
graphs of the embedding of $D$ in $F(D_s)$.

\begin{figure}[h]
\begin{center}
\includegraphics[height=1.75in]{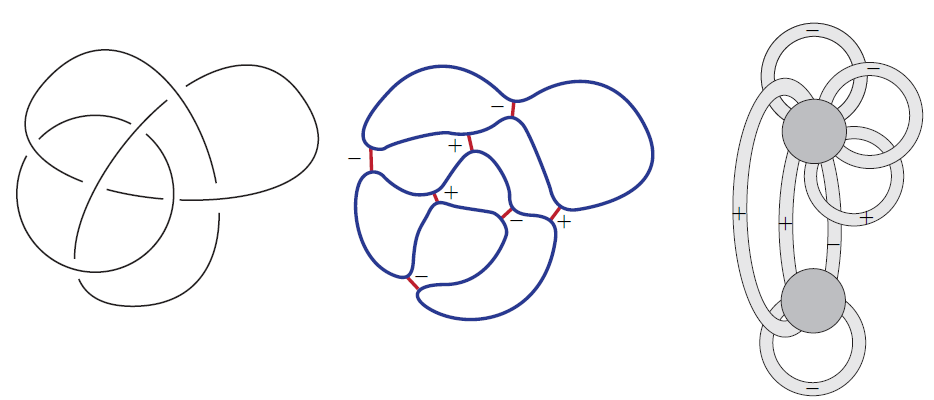}
\end{center}
\caption{Ribbon graph corresponding to any state of a link diagram (figure from \c{Moffatt_duals}).}
\label{figure:ribbon-graph-from-state}
\end{figure}

We can also obtain $G_s$ (and $G_{\overline{s}}$) directly from the
link diagram as described above. We add a $+$ sign for the edges
obtained by a $B$-smoothing and a $-$ sign for the edges obtained by
an $A$-smoothing. This gives us a signed ribbon graph which keeps
track of the smothings.  Note that the sign convention is chosen so
that the signs on edges agree with those for the Tait graph in the
case when the state is the ``Tait'' state i.e $F(D_s)=S^2$. See Figure
\ref{figure:ribbon-graph-from-state} \cite{Moffatt_duals}.

With this construction, any link diagram $D$ with $n$ crossings gives
rise to a set of $2^n$ (signed) ribbon graphs associated to the states of $D$.
It turns out that all of these ribbon graphs are {\em partial duals}
of each other, in the sense of Chmutov \c{Chmutov}.  It follows that a
ribbon graph represents a link diagram if and only if it is a partial
dual of a plane graph.  Recently in \c{moffatt_minors}, Moffatt used
this fact to completely characterize such ribbon graphs in terms of
three excluded minors.  

\subsection{Quasi-trees}

\begin{center}
\begin{figure}[h]
\includegraphics[height=1.05 in]{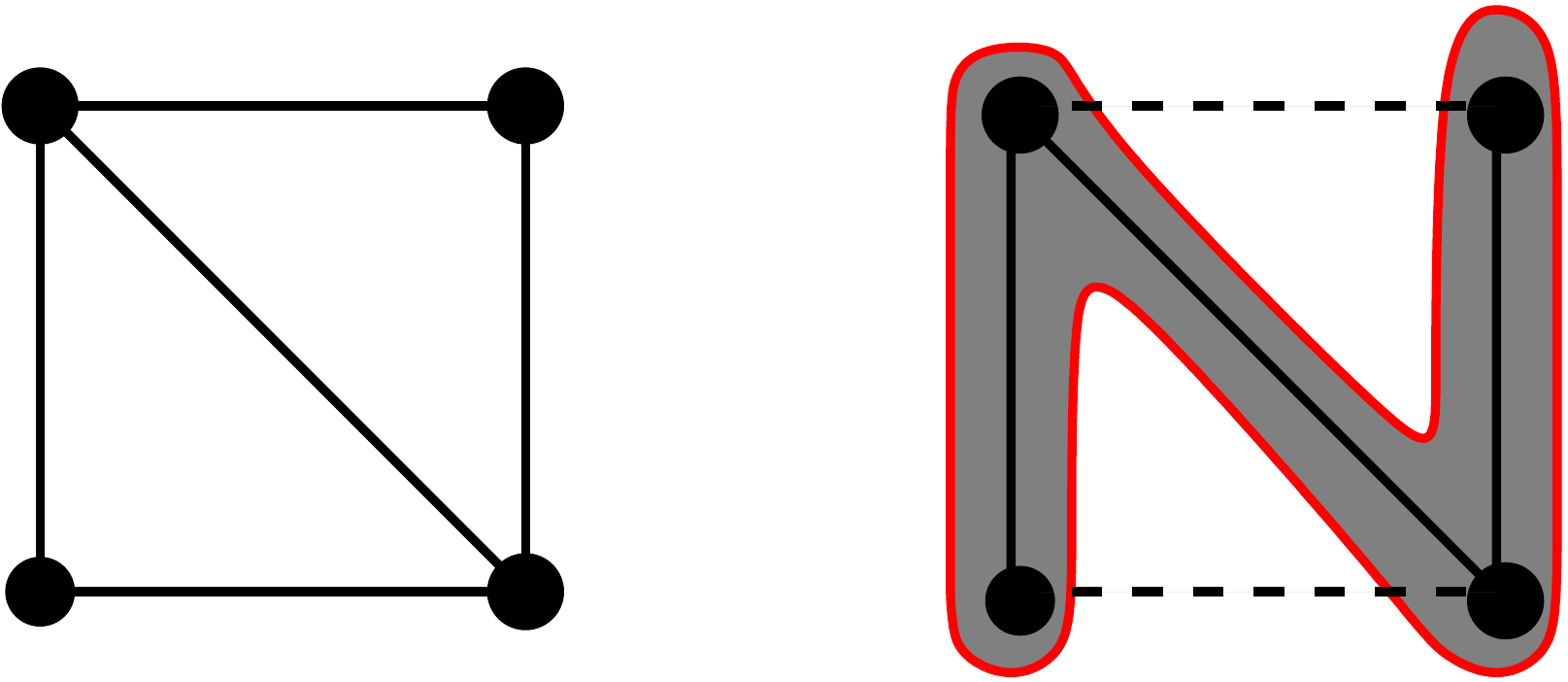} \hspace{0.75in} 
\includegraphics[height=1.2 in]{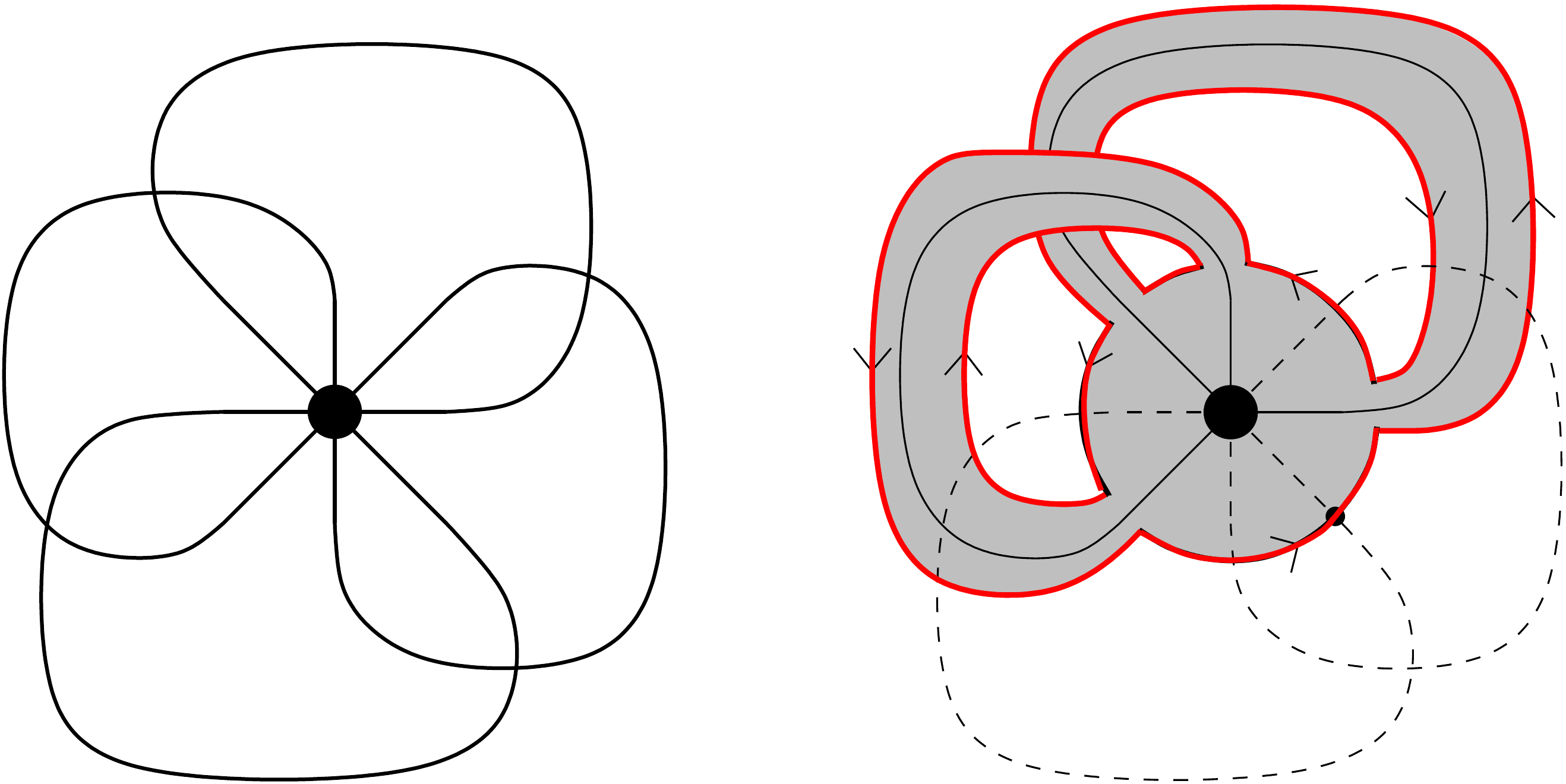}
\caption{Regular neighborhoods of spanning trees and quasi-trees}
\label{figure:st-qt}
\end{figure}
\end{center}

A {\it quasi-tree} $\Q$ of a ribbon graph
$\G$ is a ribbon spanning subgraph  with $f(\Q)=1$.
So a quasi-tree is a spanning ribbon subgraph whose regular
neighborhood on $S(\G)$ has one boundary component. This generalizes
the analogous defining property of a spanning tree of a plane graph.
See Figure \ref{figure:st-qt}.

The Tutte polynomial counts the number of spanning trees of a
connected graph $G$ by the specialization $T_G(1,1)$.  (For any
alternating knot, this is exactly the determinant of the knot.)  For
any ribbon graph, we proved with Stoltzfus that the Bollob\'as--Riordan--Tutte
polynomial also counts the number of quasi-trees of every genus by a
specialization as follows.

\begin{prop}[\cite{qtbrt}]
\label{qtcount}
Let $q(\G;t,Y)=C(\G;1,Y,t Y^{-2})$. Then $q(\G;t,Y)$ is a polynomial
in $t$ and $Y$ such that 
$$ q(\G;t) := q(\G;t,0)=\sum\nolimits_j a_j t^j$$ 
where $a_j$ is the number of quasi-trees of genus $j$. 
Consequently, $q(\G;1)$ equals the number of quasi-trees of $\G$. 
\end{prop}

\begin{question}\label{qpoly}
Let $\G$ be the all--$A$ ribbon graph for a diagram $D$ of a knot $K$.  If $g_T(D)=g_T(K)$, is $q(\G;t)$ an invariant of $K$?
\end{question}

\section{Turaev genus and knot homology} \label{KHsection}

In his theorem mentioned above, Thistlethwaite \c{Thistlethwaite} gave
an expansion of the Jones polynomial $V_L(t)$ in terms of spanning
trees of any Tait graph $G$ of $L$. 
In \c{KH}, for any connected link diagram $D$, we defined the spanning
tree complex $\C(D)=\{\C_v^u(D), \del\}$, whose generators correspond
to spanning trees $T$ of $G$, and whose homology is the reduced
Khovanov homology.  As described precisely below, $\C(D)$ is at most
$g_T(L)+1$--thick, where $g_T(L)$ is the Turaev genus.

\subsection{Spanning tree expansion} 
We first describe the spanning tree expansion for the Jones
polynomial, and then the spanning tree expansion for Khovanov
homology, which is also similar to the one for knot Floer homology.

Fix an order on the edges of $G$.  For every spanning tree $T$ of $G$,
each edge $e\in G$ has an activity with respect to $T$, as follows.
If $e \in T$, $\mathit{cut(T,e)}$ is the set of edges that connect
$T\setminus e$.  If $f \notin T$, $\mathit{cyc(T,f)}$ is the set of
edges in the unique cycle of $T \cup f$.  Note $f \in cut(T,e)$ if and
only if $e \in cyc(T,f)$.  An edge $e \in T$ (resp. $e \notin T$) is
{\em live} if it is the lowest edge in its cut (resp. cycle), and
otherwise it is {\em dead}.

For any spanning tree $T$ of $G$, the {\em activity word} $W(T)$ gives
the activity of each edge of $G$ with respect to $T$.  The letters of
$W(T)$ are as follows: $L,\ D,\ \ell,\ d$ denote a positive edge that
is live in $T$, dead in $T$, live in $G-T$, dead in $G-T$,
respectively; $\bar{L},\ \bar{D},\ \bar{\ell},\ \bar{d}$ denote
activities for a negative edge.  Note that $T$ is given by the capital letters of $W(T)$.

Thistlethwaite assigned a monomial $\mu(T)$ to each $T$ as follows:
$$ L^p D^q \ell^r d^s \bar{L}^x \bar{D}^y \bar{\ell}^z \bar{d}^w \quad \Rightarrow \quad \mu(T)= (-1)^{p+r+x+z}A^{-3p+q+3r-s+3x-y-3z+w} $$
\begin{theorem}[\c{Thistlethwaite}]\label{ThistThm}
  Let $G$ be the Tait graph of any connected link diagram $D$ with any
  order on its edges.  Let $\kb{D}$ denote the Kauffman bracket
  polynomial of $D$. Summing over all spanning trees $T$ of $G$,
  $\kb{D}= \sum_T \mu(T)$.
\end{theorem}

The activity word $W(T)$ contains much more information than just
$\mu(T)$.  A \emph{twisted unknot} $U$ is a diagram of the unknot
obtained from the round unknot using only Reidemeister I moves.
$W(T)$ determines a twisted unknot $U(T)$ by changing the crossings of
$D$ according to Table \ref{Table1} for dead edges, and leaving the
crossings unchanged for live edges (Lemma 1 \c{KH}).  In Table
\ref{Table1}, the sign of the crossing in $U(T)$ is indicated for
unsmoothed crossings, and Kauffman state markers are indicated for
smoothed crossings.
\begin{table}[h]
\caption{Activity word for a spanning tree determines a twisted unknot}
\label{Table1}
\begin{center}
\begin{tabular}{cc|cc|cc|cc}
$L$ & $D$ & $\ell$ & $d$ & $\bar{L}$ & $\bar{D}$ &
$\bar{\ell}$ & $\bar{d}$ \\
\hline
$-$ & $A$ & $+$ & $B$ & $+$ & $B$ & $-$ & $A$ \\
\includegraphics[height=0.5cm]{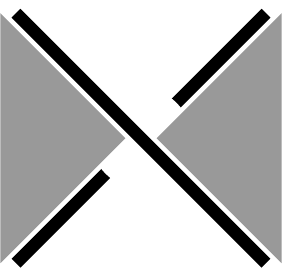} & \includegraphics[height=0.5cm]{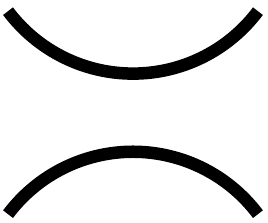} &
\includegraphics[height=0.5cm]{poscross} & \includegraphics[height=0.5cm]{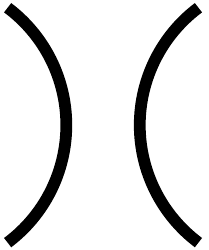} &
\includegraphics[height=0.5cm]{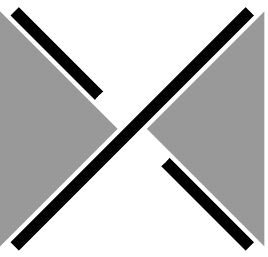} & \includegraphics[height=0.5cm]{posA} &
\includegraphics[height=0.5cm]{negcross} & \includegraphics[height=0.5cm]{posB} \\
\end{tabular}
\end{center}
\end{table}

\begin{figure}[h]
\includegraphics[height=1in]{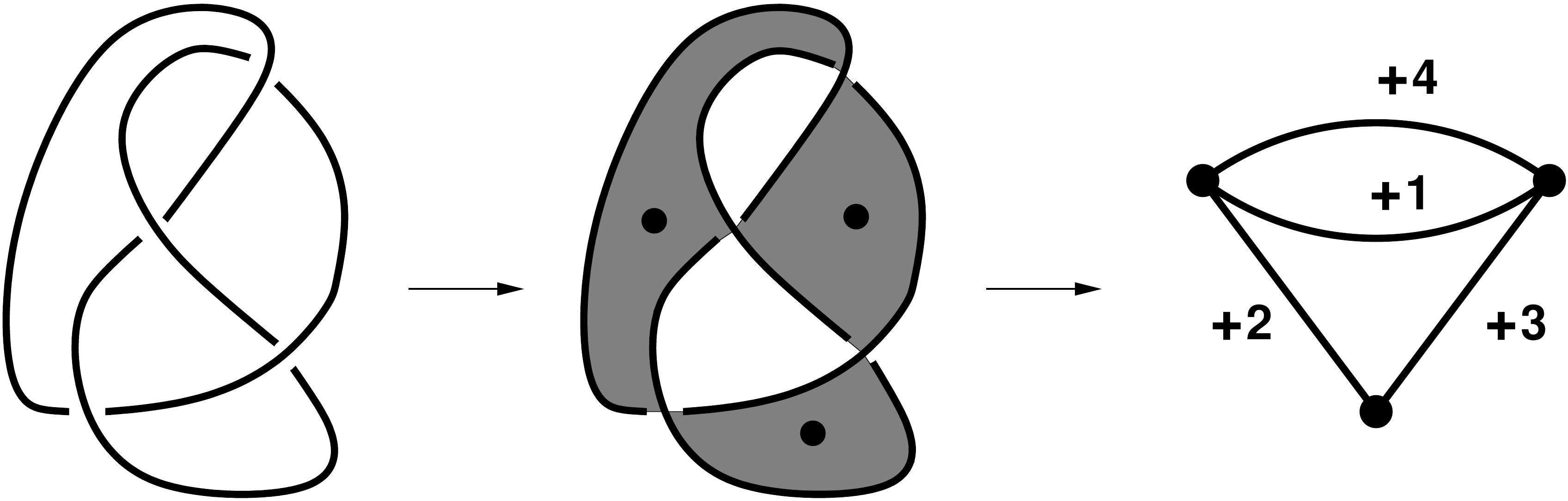}\\
\ \\
\includegraphics[height=4.5in]{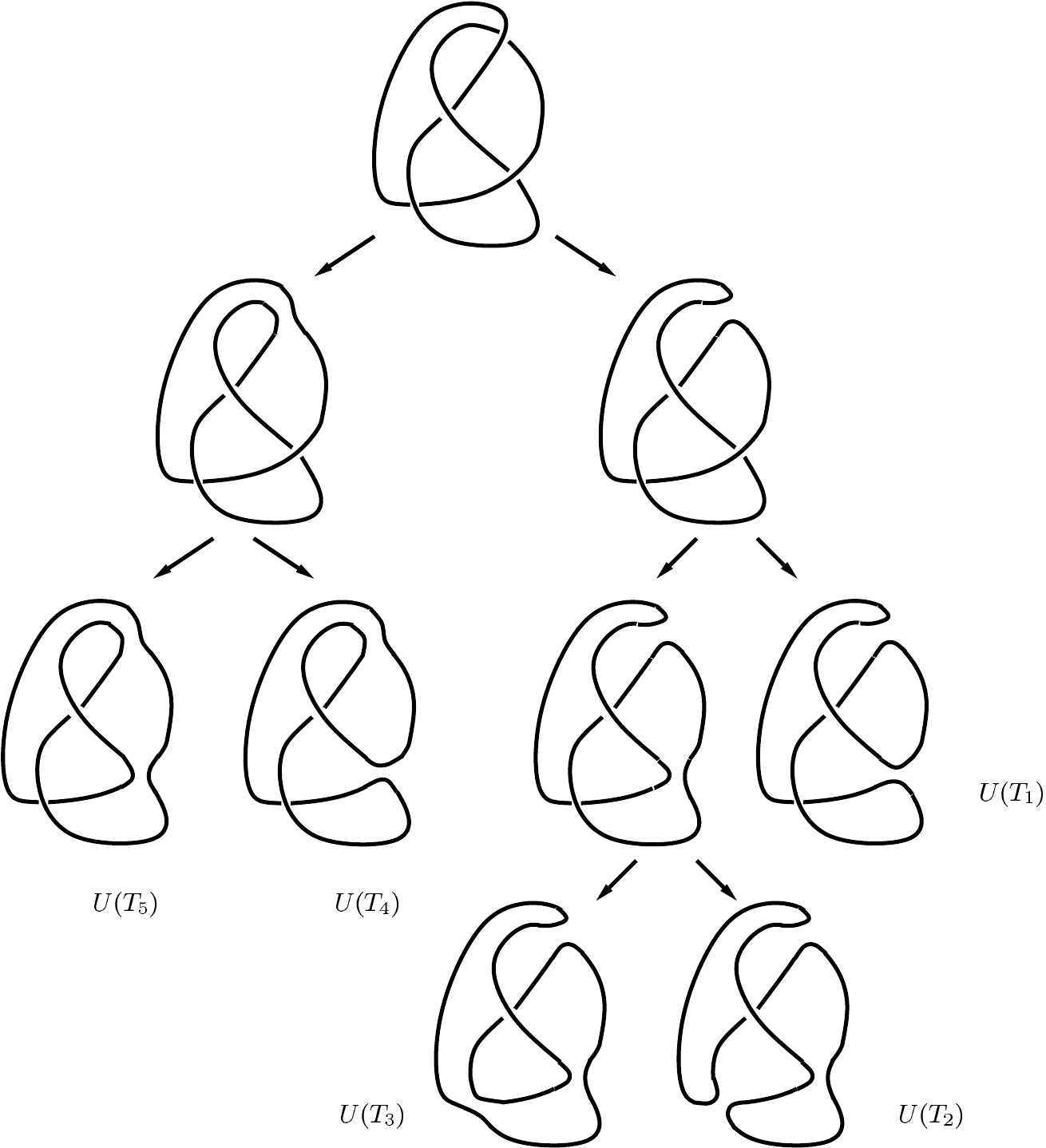}
\caption{Twisted unknots corresponding to spanning trees in Table
  \ref{table:unknots}. Crossings are smoothed in reverse order of
  the crossings of the diagram.  At every node of the skein resolution tree,
  $A$--smoothings are on the left, and $B$--smothings are on the right.}
\label{figure:unknots}
\end{figure}

We can also consider each $U(T)$ as a partial smoothing of $D$
determined by $W(T)$.  In fact, there exists a skein resolution tree
for $D$ whose leaves are exactly all the partial resolutions $U(T)$,
for each spanning tree $T$ of $G$.  Let $\sigma(U) = \# A$-smoothings
$- \#B$-smoothings, and let $w(U)$ be the writhe.  If $U$ corresponds
to $T$, then $\mu(T)=A^{\sigma(U)}(-A)^{3w(U)}$ is exactly the
monomial above Theorem \ref{ThistThm}.  As Louis Kauffman pointed out,
this is how humans would compute $\kb{D}$: Instead of smoothing all
the way to the final Kauffman states, a human would stop upon reaching
any twisted unknot $U$, and use the formula $\mu(T)$.  We illustrate
all of this for the figure-eight knot diagram in Figure 
\ref{figure:unknots} and Table \ref{table:unknots}.

\begin{table}[h]

\begin{tabular}{|c|c|c|c|c|c|}
\hline
&&&&& \\
Spanning trees  & 
\includegraphics[height=0.4in]{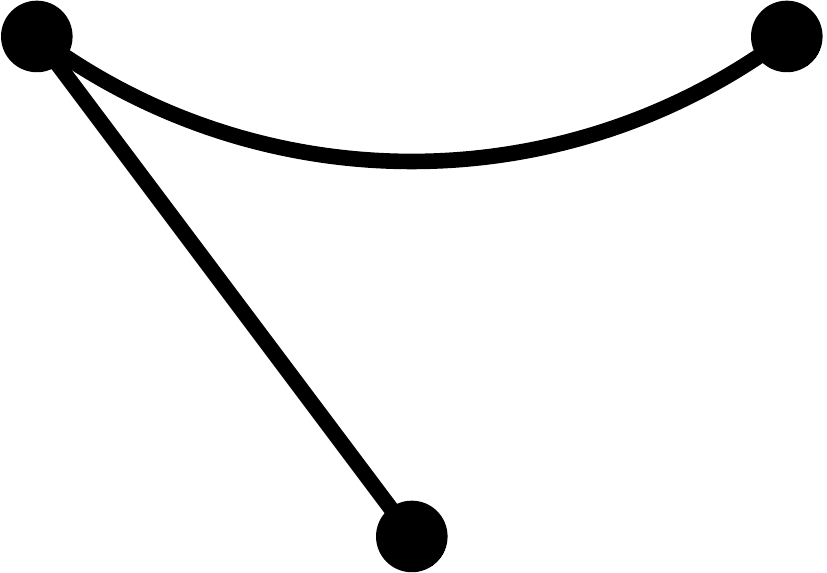} &
\includegraphics[height=0.4in]{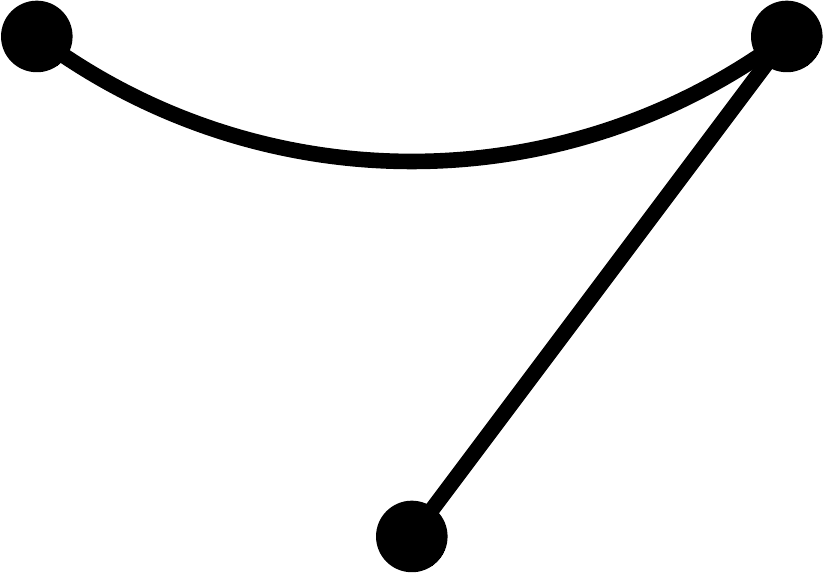} &
\includegraphics[height=0.4in]{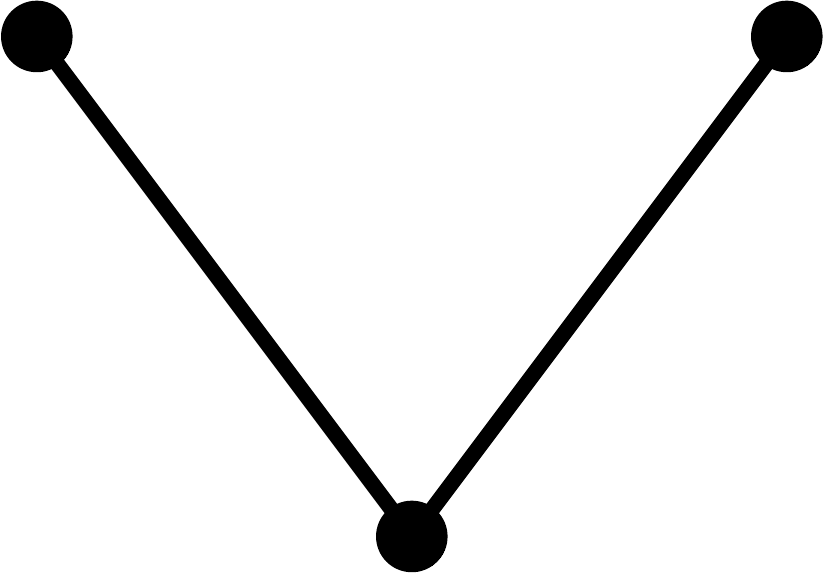} &
\includegraphics[height=0.4in]{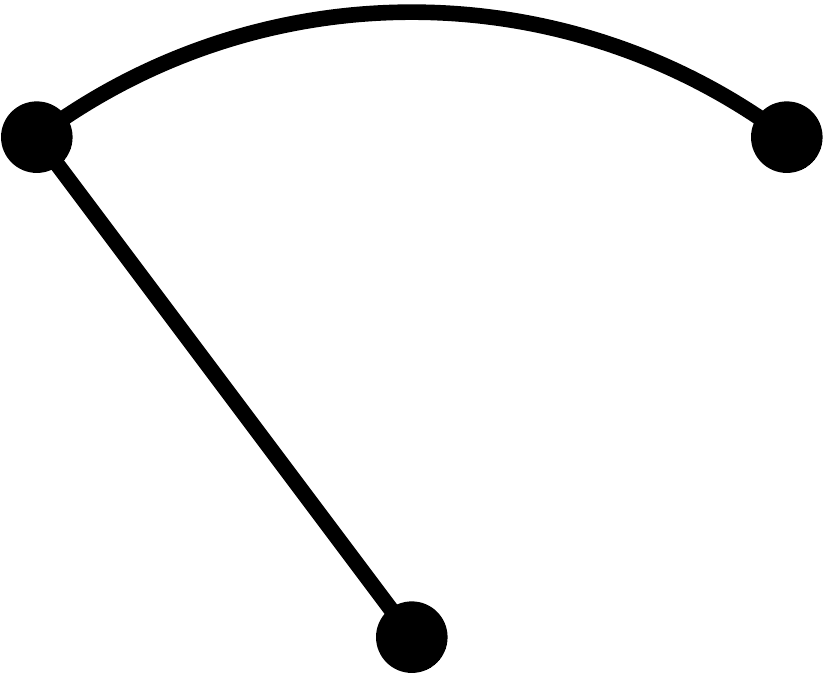} &
\includegraphics[height=0.4in]{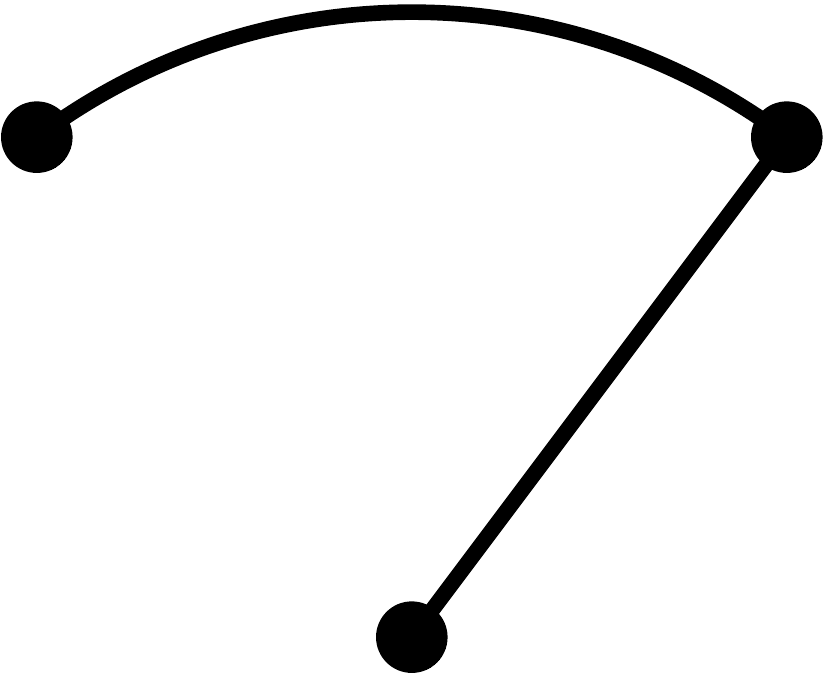} \\
&&&&& \\

& $T_1$ & $T_2$ &$T_3$ &$T_4$ &$T_5$ \\  
&&&&& \\
\hline
&&&&& \\
Activities & $LLdd$ & $LdDd$ &$\ell DDd$ &$\ell LdD$ &$\ell \ell DD$ \\  
&&&&& \\
\hline
&&&&& \\
Weights & $A^{-8}$ & $-A^{-4}$ &$-A^4$ &$1$ &$A^8$ \\  
&&&&& \\
\hline
\end{tabular}
\caption{Spanning trees of the Tait graph of the figure-8 knot. They correspond 
to twisted unknots in Figure \ref{figure:unknots}, and
their weights add up to $\langle D \rangle=A^{-8}-A^{-4}+1-A^4+A^8$.}
\label{table:unknots}

\end{table}

For any connected link diagram $D$, we choose the checkerboard
coloring such that its Tait graph $G$ has more positive edges than
negative edges, and in case of equality that the unbounded region is
unshaded.  In \c{KH}, we defined the spanning tree complex
$\C(D)=\{\C_v^u(D), \del\}$, whose generators correspond to spanning
trees $T$ of $G$.  The $u$ and $v$--grading are determined by $W(T)$
as follows:
$$ u(T)  =   \# L - \# \ell - \# \bar{L} + \# \bar{\ell} \quad {\rm and }\quad v(T) =  \# L + \# D = e_+(T) $$
\begin{theorem}[\c{KH}]\label{KHthm}
  For any connected link diagram $D$, there exists a spanning tree
  complex $\C(D)=\{\C_v^u(D), \del\}$ with $\del$ of bi-degree $(-1,-1)$ that is a deformation
  retract of the reduced Khovanov complex.
\end{theorem}

\subsection{Turaev genus and Khovanov homology}
The key idea for relating Khovanov homology to Turaev surfaces, Turaev
genus and ribbon graphs is our observation that there is a one-to-one
correspondence between spanning trees of the Tait graph and
quasi-trees of the all-A ribbon graph \c{dkh}.

Let $D$ be a connected link diagram, $G$ be its Tait graph for which
the number of positive edges is greater than or equal to the number of
negative edges and let $\G$ be the all-A ribbon graph.  In \c{dkh}, we
proved with Stoltzfus

\begin{theorem}[\c{dkh}]\label{dkh_thm}
Quasi-trees of $\G$ are in one-one correspondence with spanning trees of $G$:
$$ \Q_j \leftrightarrow \TT_v \qquad {\rm where }\qquad v+j = \frac{v(G) + e_+(G) - v(\G)}{2} $$
$\Q_j$ denotes a quasi-tree of genus $j$, $T_v$ denotes a spanning
tree with $v$ positive edges, and $e_+(G)$ equals the number of positive edges 
in $G$. 
\end{theorem}

To construct the spanning tree chain complex in \c{KH}, every spanning
tree $T$ of the Tait graph $G$ was given a bigrading $(u(T),
v(T))$.  By Theorem \ref{dkh_thm}, the $v$-grading, which is the
number of positive edges in $T$, is determined by the genus of the
corresponding quasi-tree $\Q$.  The $u$-grading, which was defined
using activities in the sense of Tutte, also has a quasi-tree analogue
in terms of the ordered chord diagram for $\Q$.

If $D$ has $n$ ordered crossings, let $\G$ be given by permutations
$(\s_0,\s_1,\s_2)$ of the set $\{1,\ldots,2n\}$, such that the $i$-th
crossing corresponds to half-edges $\{2i-1, 2i\}$, which are marked on
the components of the all--$A$ state of $D$.  
We give the components of the all--$A$
state of $D$ the admissible orientation for which outer ones are
oriented counterclockwise (see \c{DFKLS1}).  In this way, every
component has a well-defined positive direction.

The orbits of $\s_0$ form the vertex set. In particular, $\s_0$
is given by noting the half-edge marks when going in the positive
direction around the components of the all--$A$ state of $D$.
The other permutations are given by
$ \s_1 = \prod_{i=1}^{n} (2i-1, 2i)$ and $\s_2 = \s_1 \circ \s_0^{-1}$

\begin{figure}
\includegraphics[height=1.5in]{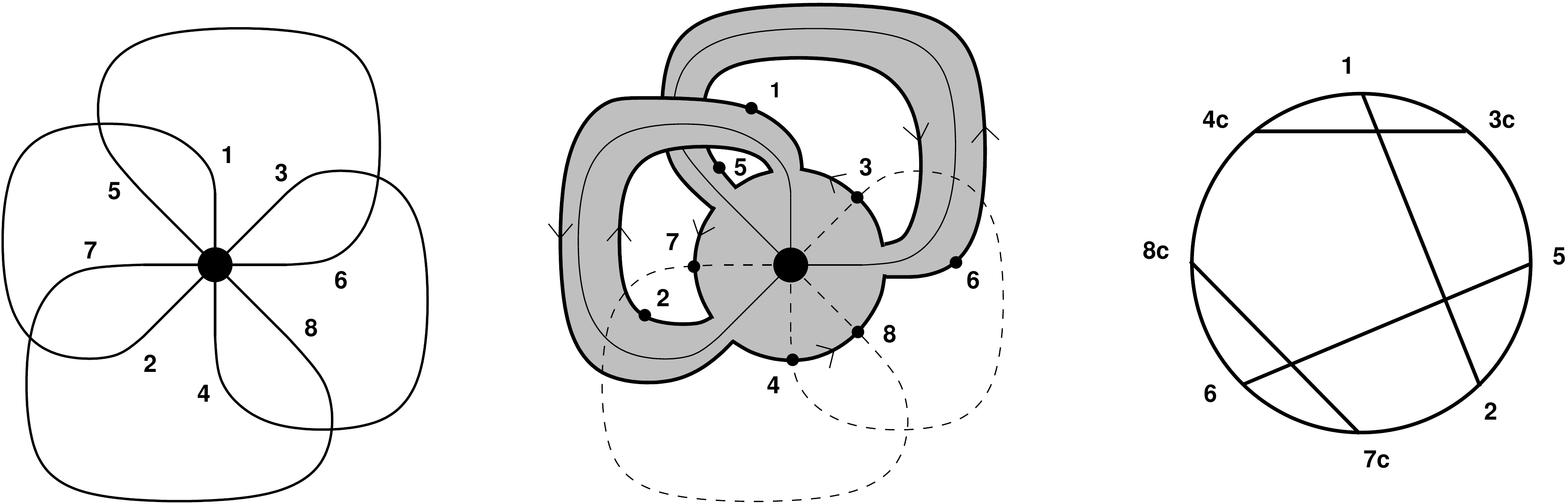}
\caption{Ribbon graph $\G$, quasi-tree $\Q=(12)(56)$ with curve 
$\gamma_{\Q}$, chord diagram $C_{\Q}$}.
\label{figure:qt_face}
\end{figure}

Let an {\em ordered chord diagram} denote a circle marked with
$\{1,\ldots,2n\}$ in some order, and chords joining all pairs $\{2i-1, 2i\}$.
By Proposition 1 of \c{dkh}, every quasi-tree $\Q$ corresponds to the ordered chord diagram
$C_{\Q}$ with consecutive markings in the positive direction given
by the permutation:
$$\s(i) =
\begin{cases}
\s_0(i) & i\notin\Q \\
\s_2^{-1}(i) & i\in\Q
\end{cases} $$
For example see Figure \ref{figure:qt_face}.

Using $\min(i,\s_1(i))$, there is an induced total order on the chords
of $C_{\Q}$. A chord is {\em live} if it does not intersect
lower-ordered chords, and otherwise it is {\em dead}.  For any
quasi-tree $\Q$, an edge $e$ is {\em live} or {\em dead} when the
corresponding chord of $C_{\Q}$ is live or dead.  In Figure
\ref{figure:qt_face}, we show $C_{\Q}$ such that the only edge live
with respect to $\Q$ is $(12)$.

To compute the genus $g(\Q)$ from $C_{\Q}$, let $C$ be the sub-chord diagram of
chords that correspond to edges in $\Q$.  Then $g(\Q)$ is half the rank
of the adjacency matrix of the intersection graph of $C$ \c{BR1}.

In \c{dkh}, we proved that if the spanning tree $T$ corresponds to
$\Q$, as in Theorem \ref{dkh_thm}, then the chord diagram $C_{\Q}$
parametrizes the regular neighborhood of $T$ formed by the
appropriate smoothings of $D$.  Consequently, we proved that the
$i$-th edge of $\G$ is live with respect to $\Q$ if and only if the
$i$-th edge of $G$ is live with respect to $T$.  This is the essential
reason that the spanning tree complex can be expressed entirely as a
bigraded quasi-tree complex, with the Turaev genus as one of the
gradings:
\begin{theorem}[\c{dkh}]  
For a knot diagram $D$ with all--$A$ ribbon graph $\G$, there exists a quasi-tree complex
$\C(\G)=\{\C_v^u(\G), \del\}$ that is a deformation retract of the
reduced Khovanov complex, where
$ \C_v^u(\G) = \Z\kb{\Q\subset \G |\; u=u(\Q),\; v=-g(\Q)}$.
\end{theorem}

\begin{corollary}[\c{dkh}]\label{KHthick}
For any knot $K$, the width of its reduced Khovanov homology
$$ w_{KH}(K) \leq 1+g_T(K).$$
\end{corollary}

The proof follows from the fact that $\displaystyle w_{KH}(K)\leq
\max_{T\subset G} v(T) - \min_{T\subset G} v(T) + 1$, where $G$ is the
Tait graph of any diagram of $K$, and that for the all--$A$ ribbon graph $\G$,
$$ g(\G)  = \max_{\Q\subset\G} g(\Q) - \min_{\Q\subset\G} g(\Q) = \max_{T\subset G} E_+(T) - \min_{T\subset G} E_+(T) = \max_{T\subset G} v(T) - \min_{T\subset G} v(T). $$

Note that $u(\Q) = -w(U(T))$, where the twisted unknot $U(T)$ comes
from the spanning tree corresponding to the quasi-tree $\Q$.  The
grading $v(\Q)$ is related to Rasmussen's $\delta$--grading for
Khovanov homology as $\delta = 2v+k$, where $k$ is a constant that
depends only on $D$.
Because the $(u,v)$ gradings are linear combinations of Khovanov's $(i,j)$
gradings, the width refers to the diagonals of Khovanov's complex, as
in the following figure from \c{DasbachLowrance2}:
\begin{center}
\includegraphics[height=2.5in]{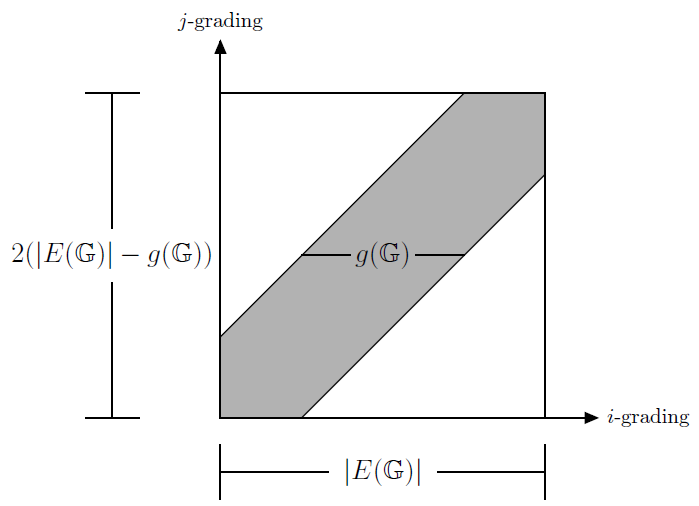}
\end{center}

The rational Khovanov homology of $(3,q)$--torus links was
computed in \c{Turner}, where it was shown that the width of the 
Khovanov homology of torus knots of type $(3,3N+1)$ and $(3,3N+2)$ is
exactly $N+2$. Corollary \ref{KHthick} gives a family 
of links with unbounded Turaev genus. 
 
\begin{corollary} [\c{dkh}]\label{tg-unbounded}
The Turaev genus of $(3,q)$--torus knots is unbounded.
\end{corollary}

Lowrance \c{Lowrance2011} and Watson \c{watson} have proved that the
width of Khovanov homology remains unchanged after replacing a
crossing in a link diagram with an alternating rational tangle,
provided the crossing satisfies certain conditions. Using 
Corollary \ref{tg-unbounded}, this generates many families of knots 
with unbounded Turaev genus.

\subsection{Turaev genus and knot Floer homology}
Finally, we turn briefly to knot Floer homology $\widehat{HFK}$ (see
the recent surveys \cite{juhasz, Manolescu}).  Lowrance \cite{lowrance1} proved
the analogous bound to Corollary \ref{KHthick}. 

\begin{theorem}[\c{lowrance1}]\label{HFthick}
    Let $K \subset S^3$ be a knot. The knot Floer width of $K$ is bounded by
    the Turaev genus of $K$ plus one: $w_{HF} (K) \leq  g_T (K) + 1$.
\end{theorem}

This bound follows from the same idea as for the Khovanov homology. In
\cite{os03}, Ozsv\'ath and Szab\'o showed that for any diagram $D$ of
a knot $K$, there exists a complex whose generators are in one-to-one
correspondence with spanning trees of the Tait graph of $D$ and whose
homology is the knot Floer homology of $K$.  The $\delta$--grading on
$\widehat{HFK}$ is defined as the difference of the Alexander and
homological gradings. The $\delta$-grading of a spanning tree when
considered in the reduced Khovanov complex is the same as the
$\delta$-grading of that spanning tree when considered in the knot
Floer complex (see \cite{DasbachLowrance1}). Thus, by Theorem
\ref{dkh_thm}, the width of the spanning tree complex giving knot
Floer homology is also bounded by $g_T(K)+1$, hence so is the width of
the homology.
 
\section{Bounds for the Turaev genus}\label{bounds}

The first bound that was discovered, equation (\ref{spaneq}) in
Section 2, motivated interest in the Turaev genus. Corollary \ref{KHthick}
and Theorem \ref{HFthick} give lower bounds for 
the Turaev genus in terms of homological width. Here, we survey other known 
bounds for the Turaev genus.

\subsection{Alternating embeddings of link diagrams} 

Let $g_F(L)$ be the minimal genus of an unknotted closed orientable
surface $F$ on which $L$ can be cellularly embedded on $F$, i.e. the
complementary regions $F-L$ are discs, and such that $L$ has an
alternating diagram on $F$.  We saw in Section 1 that any
non-alternating pretzel link $P$ can be made alternating on the torus,
so $g_F(P)=1$.  Colin Adams called such links ``toroidally
alternating.''  Using Turaev's construction it's easy to see that
$g_T(P)=1$ for such pretzel links. (See \cite{Hayashi} for properties
of cellularly embedded links with alternating diagrams on higher genus
surfaces.)  In general, $$ g_F(L) \leq g_T(L) $$

There are examples due to Adam Lowrance of a family of links $L_n$ for
which $g_F(L_n)=1$ but $g_T(L_n)\to\infty$ as $n\to\infty$.  The idea
is to start with torus links $T_n$ for which $w_{KH}(T_n)\to\infty$, and
then to insert a small alternating tangle to get $L_n$ so that $g_F(L_n)=1$ but
$w_{KH}(T_n) = w_{KH}(L_n)$. 

Note that these provide examples of cellularly embedded alternating links
on an unknotted closed orientable surface, but they do not satisfy the
Morse decomposition property of the Turaev surface, which is property
(e) listed in Section 1.

\subsection{Dealternating number} Let $dalt(L)$ be the {\em
  dealternating number}, which is defined to be the minimal number of
crossing changes needed to make a diagram of $L$ into an alternating
diagram.  When $dalt(L)=1$, the link is called {\em almost
  alternating}.  Abe and Kishimoto \cite{AbeKishimoto} proved $$
g_T(L) \leq dalt(L). $$ Moreover, the almost alternating torus knots
are exactly the only two knots which are both torus and pretzel knots,
$T(3,4)=P(3,3,-2)$ and $T(3,5)=P(5,3,-2)$.

\subsection{Homological invariants}
For a knot $K$, let $\sigma(K),\ s(K)$, and $\tau(K)$ denote the
signature of $K$, the Rassmussen $s$--invariant of $K$ which comes
from Khovanov homology, and the $\tau$--invariant of $K$ which comes
from the knot Floer homology.  If $K$ is any alternating knot, then 
$2\tau(K) = s(K) = -\sigma(K)$.

In \c{DasbachLowrance1}, Dasbach and Lowrance proved the same lower bounds for the Turaev genus that Abe had proved for the dealternating number:
\begin{equation}\label{DLeq}
 \left|\tau(K)+\frac{\sigma(K)}{2}\right|\leq g_T(K),\quad \frac{|s(K)+\sigma(K)|}{2}\leq g_T(K),\quad \left|\tau(K)-\frac{s(K)}{2}\right|\leq g_T(K) 
\end{equation}
These are known to be equalities for alternating knots (when $g_T(K)=0$).  
Dasbach and Lowrance gave certain examples of $(3,q)$--torus knots for which they are not sharp.

For any link $L$, let $L^{(r)}$ be its $r$--fold parallel.  Huggett, Moffatt and Virdee \c{HMV} gave an upper bound on the Turaev
genus of $L^{(r)}$, 
$$ g_T(L^{(r)})\leq (r + 1)\cdot g_T(L) + r^2c - r $$
where $c$ is the crossing number of any diagram $D$ for which $g_T(D)=g_T(L)$.

\section{Open questions and research directions}\label{open}

\subsection{Adequate knots}

Equations (\ref{spaneq}) and (\ref{adeqeq}) in Section 2, focused
attention on adequate diagrams because it follows that they have
minimal crossing number.  These are much more general than alternating
diagrams.  For example, the $r$--fold parallel of an adequate diagram
is adequate.

Recall that we say that a diagram $D$ is $A$--adequate if $|s_A|>|s|$
for any state $s$ with exactly one $B$--smoothing, and is
$B$--adequate if $|s_B|>|s|$ for any state $s$ with exactly one
$A$--smoothing.  In terms of ribbon graphs, $D$ is $A$--adequate
($B$--adequate) if $\G_A\ (\G_B)$ has no loops.  $D$ is adequate if it
is both $A$--adequate and $B$--adequate.  A knot or link is adequate
if it has an adequate diagram.

For an adequate knot, Abe \c{Abe1} proved that the inequality in
Corollary \ref{KHthick} is an equality.  Thus, extending equation
(\ref{adeqeq}) for an adequate knot $K$,
$$ g_T(K) = w_{KH}(K) - 1 = c(K) - {\rm span}\, V_K(t).$$ 

\begin{question}
Is Lowrance's analogous inequality for knot Floer homology an equality for adequate knots? 
\end{question}
\begin{question}
For any two knots $K$ and $K'$, is $g_T(K\#K')=g_T(K)+g_T(K')$?
\end{question}
\begin{question}
If $K$ and $K'$ are mutant knots, is $g_T(K)=g_T(K')$?  
\end{question}
In \c{Abe1}, Abe answered both questions for adequate knots, proving that the Turaev genus is additive
under connect sum, and is invariant under mutation for adequate knots.  (Mutation of an adequate diagram preserves adequacy.)

It is not known whether the other inequalities above are equalities for adequate knots:
\begin{question}
If $K$ is an adequate knot, is $g_T(K) = g_F(K)$?
\end{question}
\begin{question}
Are the inequalities (\ref{DLeq}) equalities for adequate knots?
\end{question}

In \cite{ck09}, we considered an operation on diagrams to extend a
twist on two strands by any rational tangle, as defined by Conway.
This operation can change the link type, but it preserves the
properties of the diagram $D$ being alternating, adequate, or
quasi-alternating \cite{ck09}, and it also preserves $g_F(D)$
and $g_T(D)$.  

\begin{question}
Let $L_1$ and $L_2$ be links whose diagrams are related by extending a twist on two strands by some rational tangle.
Is $g_F(L_1)=g_F(L_2)$ and is $g_T(L_1)=g_T(L_2)$?
\end{question}

We now turn to a related open problem in knot homology.  Because knot
Floer homology detects the Seifert genus, it can detect mutation; the
Seifert genus of the Conway knot is 3, and that of the
Kinoshita-Teresaka knot is 2.  It is not known whether Khovanov
homology is invariant under mutation (odd Khovanov homology is known
to be invariant).  For both homology theories, though, the rank of the
homology in each $\delta$--grading (i.e., the $v$--grading discussed
above for Khovanov homology) is conjectured to be invariant under
mutation (see Conjecture 3 of \cite{juhasz}).  Following Abe's
results, and by the proofs of Corollary \ref{KHthick} and the similar
result by Lowrance for knot Floer homology, we are led naturally to
the following conjecture.  For an adequate knot, it would imply the
mutation invariance of the ranks of both homology theories in each
$\delta$--grading.

\begin{conjecture}
Let $K$ be an adequate knot.
\begin{enumerate}
\item For any adequate diagram of $K$, the ranks of both Khovanov and
  knot Floer homology in each $\delta$--grading is given by $q(\G;t)$,
  as in Proposition \ref{qtcount}.
\item If $K'$ is any mutant of $K$, then for any adequate diagrams of $K$ and $K'$, \\ $q(\G;t) = q(\G';t)$.
\end{enumerate}
\end{conjecture}

\subsection{Quasi-alternating links} Quasi-alternating links were
first defined in \c{Ozsvath-Szabo:DoubleCovers} and it was shown in
\c{Ozsvath-Manolescu, Ozsvath-Szabo:DoubleCovers} that, like
alternating links, they are homologically thin with respect to both
Khovanov and knot Floer homology.  As a result, the homological bounds
on the Turaev genus discussed above vanish for quasi-alternating
links. Examples of quasi-alternating links of Turaev genus one include
non-alternating pretzel links and, more generally, non-alternating
Montesinos links (see \c{ck09}).
To answer the following question requires a new kind of lower bound for the
Turaev genus.

\begin{question}\label{QAquestion}
  For any $g>1$, do there exist quasi-alternating links with 
Turaev genus equal to $g$?
\end{question}

\subsection{Geometry of knot complements}
For any Kauffman state $s$ of a knot $K$, the state surface $F_s$ is
constructed like a Seifert surface: state circles bound disjoint
disks, which are connected by half-twisted bands such that $\partial
F_s = K$.  The ribbon graph $\G_s$ embeds as a spine of the surface
$F_s$.  Let $F_A$ ($F_B$) denote the all--$A$ (all--$B$) state
surface.  Ozawa \cite{ozawa} proved that if $D$ is $A$--adequate
($B$--adequate), then $F_A$ ($F_B$) is essential in $S^3-K$.  (Ozawa's
result actually holds for more general diagrams, which have a state
that is both adequate and homogeneous.)  Ozawa's theorem opens the
door to geometric results.  Futer, Kalfagianni and Purcell \cite{FKP}
related certain stable coefficients of colored Jones polynomials to
fibering data and hyperbolic volume bounds using essential state
surfaces. 

Other than this connection, very little is known about the geometry of
the knot complement and the Turaev genus.  An exciting direction to
explore may be the extent to which the Turaev genus measures how the
geometry differs from that of alternating knots.
%% This would be an exciting direction to explore.

\begin{question}
Given a knot diagram, does the Turaev genus provide any additional constraint on the geometry of the knot complement?
\end{question}

\subsection*{Acknowledgments}
We would like to thank the organizers of the Quantum Topology and
Hyperbolic Geometry Conference (Nha Trang, Vietnam, May 13-17, 2013)
for their extraordinary hospitality.  We gratefully acknowledge
support by the Simons Foundation and PSC-CUNY.

\bibliographystyle{amsplain}
\bibliography{dkh}

\end{document}